\documentclass{ifacconf}

\usepackage{graphicx}      
\usepackage{natbib}        

\usepackage{amsmath, amssymb} 
\usepackage{graphicx, color}
\usepackage{url}

\newcommand{\llim}[3]{\ensuremath{ \lim\limits_{ #1 \rightarrow #2} #3 }}
\newcommand{\fclass}[2]{\ensuremath{  \mathbb{#1}^{\, #2} }}

\newcommand{\HF}[3]{\ensuremath{
		\left[ #3 \middle|
		\begin{array}{l}
			#1      \\
			#2      \\
		\end{array}
		\right]    }}

\begin{document}
\begin{frontmatter}

\title{A space-fractional reaction-diffusion system with  cylindrical symmetry} 


\author[First]{Dimiter Prodanov}

\address[First]{Lab. of Neurotechnology, PAML-LN, IICT, Bulgarian Academy of Sciences, 1431 Sofia, Bulgaria (e-mail: dimiter.prodanov@iict.bas.bg).}

\begin{abstract}                
Diffusion within porous media, such as biological tissues, exhibits departures from conventional Fick's laws, which could result in space-fractional diffusion. The paper considers a reaction-diffusion system with two spatial compartments --  a proximal one of finite radius having a source, and an outer one extending to infinity where the source is not present but first-order decay of the diffusing species takes place. The system models the foreign body reaction around an implanted electrode.  Microscopic heterogeneity inside the tissue was modeled by a space-fractional Riesz Laplacian acting on the concentration. This allows for a flexible approach when estimating transport parameters from experimental data.  The steady-state of the system is solved in terms of Hankel and Mellin transforms, resulting in a Fox H-function.  In the integer-order case, the analytical solution reduces to a superposition of modified Bessel functions of the first and second kinds.
Solutions are exhibited by numerical quadrature of the involved Bessel function integrals.

\end{abstract}

\begin{keyword}
Hankel transform, Mellin transform, Riesz Laplacian, Fox H-function, Bessel function
\end{keyword}

\end{frontmatter}

\section{Introduction}\label{sec:intro}
 
The diffusion in porous media, such as biological tissues, is characterized by deviations from the usual Fick's diffusion laws (see discussion in \cite{Postnikov2022,Metzler2020}). 
One such type of anomalous diffusion is the space-fractional one that arises as an appropriate limit of the Continuous Time Random Walk (CTRW), when jump lengths follow a L\'evy stable distribution  \citep{Metzler2000}.
The source-free diffusion equation can be derived from a non-local generalization of the flux
\[
j =  D \nabla^\beta c
\] 
where $D$ is the generalized diffusion constant,
and $c$ is the concentration of the diffusing species, by using a co-ordinate independent fractional gradient operator $\nabla^\beta$ \citep{Silhavy2019}.
Microscopic heterogeneity inside the tissue can be lumped into the fractional order parameter $\beta$.
This potentially allows for a more flexible approach when estimating transport parameters from experimental data.  
On unbounded domain, mass conservation results in the continuity equation
\begin{equation}
 \partial_t c = \nabla \cdot j = - D (-\Delta)^{\alpha} c, \quad \alpha =\frac{1+ \beta}{2}
\end{equation}
 where the symbol $(-\Delta)^{\alpha} $ denotes the Riesz fractional Laplacian operator of order $\alpha$ \citep{riesz1938}. 
 This formulation informs the physical interpretation of the space-fractional diffusion equation.  
 
 \section{The reaction-diffusion system}\label{sec:deriv}
 A one-dimensional reaction-diffusion application, modeling the spatial distribution of reacting species, was introduced in \cite{prodanov2016c},  while the numerical aspects of the two dimensional case have been investigated in \cite{Prodanov2022}. 
The present paper considers a first-order reaction-diffusion system in two dimensions
 \begin{equation}\label{eq:syst}
 \partial_t c = - D (-\Delta)^{\alpha} c + \sigma - q \, c, \quad 0<\alpha \leq 1
 \end{equation}
where
 $\sigma$ is the intensity of the spatially-extended source and 
  $q$ is an elimination rate.
 The system can support a non-trivial steady state, which will be studied in the present contribution. 
 
\begin{figure}[h!]
	\centering
	\includegraphics[width=0.6\linewidth]{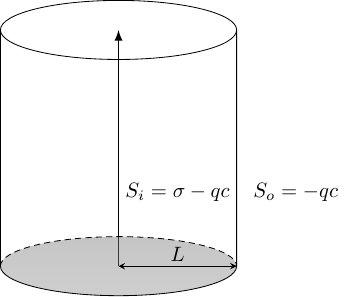}
	\caption[Model]{Geometry of the reaction-diffusion system}\label{fig:cylinder}
\end{figure}
 
We will assume cylindrical geometry of infinite extent, which effectively renders the problem 2 dimensional.
The spatial arrangement is represented in Fig. \ref{fig:cylinder}.
The system is compartmentalized into two spatial compartments -- a proximal one of   radius \textit{L}, having a source of intensity $\sigma$; 
and an outer one, which extends to infinity and where the source is not present but there is  only a first-order decay. The source term \textit{S}, therefore splits into two equations $S_i=\sigma-q c$ in the proximal compartment and $S_o= - qc$ in the distal one. 

From now on we denote the steady state concentration with the same label.
The steady state of the system  (eq. \ref{eq:syst})  can be written as
 \begin{equation}\label{eq:full2}
 	- D (-\Delta)^{\alpha} c  + \sigma (r)- q c   =0 
 \end{equation}
with boundary condition $c\left(\infty \right) = 0$. 
 For simplicity of the presentation we assume $D=1$. 
 If $D \neq 1$ then the equation can be reparametrized  as $\sigma^\prime= \sigma/D$, $q^\prime= q/D$.
  
  \section{Integer-order case}
 Since the integer-order Laplacian is local the system can split into two disjoint compartments, which interact only at their mutual boundary. 
 The full solution can be written as
 \begin{equation}
 	c (r) = \mathbf{1} (L-r) c_p (r) +  \mathbf{1} (r-L) c_d (r)
 \end{equation}
 where $\mathbf{1}()$ denotes the unit step function under the convention $\mathbf{1}(0)=1/2$.

 The equation for the proximal compartment becomes
 \begin{equation}
 - (- \Delta )^{\alpha} c_p   + \sigma   - q c_p  = 0  
 \end{equation}
 while for the outer one the original system becomes the eigen system
 \begin{equation}
  - (- \Delta )^{\alpha} c_d     - q c_d  = 0  
 \end{equation}
 In the present context, $L$ is the radius of the source, having dimension of \textit{length}.

 The solution for the integer-order case $\alpha=1$  can be obtained directly by expressing the Laplacian in cylindrical coordinates.
  One obtains the ordinary differential equation  
 \begin{equation}
  \frac{1}{r}\left(   c_p^\prime +  r c_p^{\prime\prime}  \right)   -q\, c_p = -\sigma
 \end{equation}
  where the prime denotes partial derivation in $r$.
 The left-hand side of the equation can be recognized as the modified Bessel equation of order 0.
 Therefore, for the proximal compartment it holds
 \begin{equation}
 c_p (r) = \sigma/q + k_1 I_0 (\sqrt{q} r)
 \end{equation}
 where the coefficient $k_1$ is determined from the boundary condition at $r=L$. 
 
\begin{figure}[h]
	\centering
	
	\includegraphics[width=1.0\linewidth]{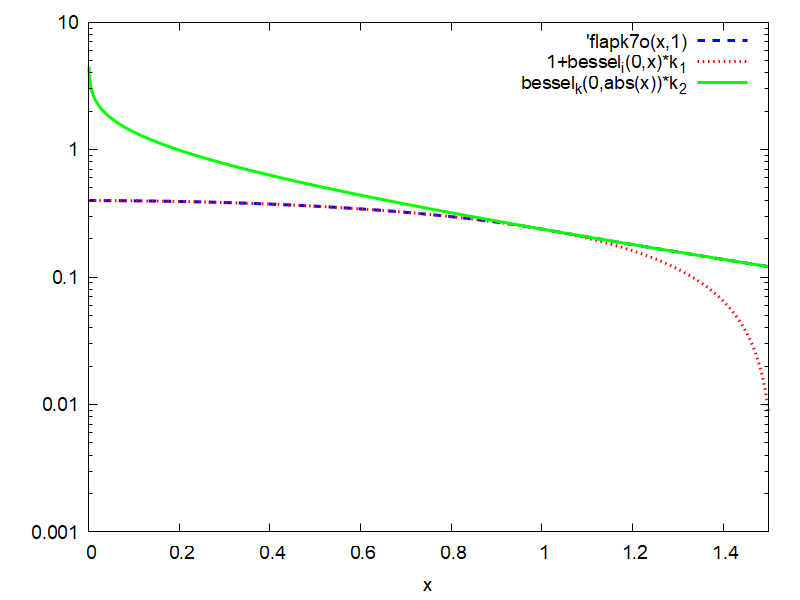}  
	
	The full solution (eq. \ref{eq:full}) is compared to the inner   (eq. \ref{eq:prox}) and  outer  (eq. \ref{eq:kf}) solutions for $\alpha = 1$.
	
	\caption{Comparison of the integer-order solutions}
	\label{fig:flapk3}
\end{figure}

The solution for the outer compartment then simply is
 \begin{equation}
 c_d=   k_2 K_0 (\sqrt{q} r)
 \end{equation}
 for the  indeterminate coefficient $k_2$.
 In the above equations, $I_0$ and $K_0$ denote the modified Bessel functions of order 0. 
 To ensure smoothness over the boundary between compartments at $r=L$ we impose matching conditions:
 \[ 
 	c_p(L) = c_d (L),  \quad
 	c^\prime_p(L) = c^\prime_d (L) 
\]
 This results in a linear system for the  coefficients, which can be solved as
  \begin{flalign}
 k_1 &=-\frac{K_1\left( L \sqrt{q}\right)  \sigma}{ P q},\\
k_2  & =
 \frac{I_1\left( L \sqrt{q}\right)  \sigma}{P q} 
  \end{flalign}
  where
  \begin{multline}
  	 P = K_0\left( L \sqrt{q}\right)  I_1\left(  L \sqrt{q}\right) +I_0\left( L \sqrt{q}\right) K_1\left(  L \sqrt{q}\right) = \\ \frac{1}{L \sqrt{q}}
  \end{multline}
  using the Wronskian identity. Therefore, the final form of the integer-order solution becomes
  \begin{flalign}
  	 c_p (r) &= \frac{\sigma}{q} - \frac{\sigma L }{\sqrt{q} } K_1( \sqrt{q} L) I_0 (\sqrt{q} r) \label{eq:prox} \\
  	 c_d (r) &= \quad\quad \frac{\sigma L }{\sqrt{q} } I_1( \sqrt{q} L) K_0 (\sqrt{q} r) \label{eq:kf}
  \end{flalign}
The solution is exhibited in Fig. \ref{fig:flapk3}.

\section{Fractional-order case}
For the fractional-order case, the compartmentalized approach is not directly applicable but the solution can be obtained in the Fourier domain. Denoting $|k|= \rho$, the stationary system can be transformed as
\begin{equation}
- \rho^{2\alpha } \hat{c} +   \sigma L \frac{ J_1 (\rho L)}{\rho} - q  \hat{c} = 0  
\end{equation}
from where one obtains the solution in transcendental-algebraic form
\[
\hat{c} =\sigma L  \frac{ J_1\left( \rho L \right) }{ \rho \left( \rho^{2\alpha} + q\right)  } 
\]
Since we assume axially-symmetric geometry, the solution in the spatial domain can be obtained by a  Hankel transform \citep{poularikas2010}:
\begin{equation}\label{eq:full}
	c (r) =  \sigma L  \int_{0}^{\infty}  \frac{J_0( \rho \,   r) J_1\left( \rho L \right)  }{\rho^{ 2 \alpha } +q } d \rho 
\end{equation}
From formal perspective this completes the analysis of the system.
However, the numerical inversion of the Hankel transform is a demanding problem. Therefore, we will look for analytical approximations of the integral.

\subsection{An asymptotic solution in terms of H-functions}\label{sec:mellint}
Remarkably, in the integer-order case the solution for the outer compartment can  also be obtained in a different way.
Notably, we assume a fictitious delta source of intensity $\sigma^\prime$ on the boundary $r=L$.
For such an ring-like source the stationary equation reads
 \begin{equation}
  \Delta \,  c_d - q c_d + \sigma^\prime \frac{\delta(r-L)}{r} = 0   
 \end{equation}
 Therefore, in the Hankel domain one obtains
 \[
 - \rho^{2  }\hat{c}  -q \hat{c} = -\sigma^\prime J_0(  \rho L)  \Longrightarrow \hat{c}= \frac{\sigma^\prime J_0(  \rho L) }{\rho^{2 } + q}
 \]
 Therefore,   
 \begin{equation}\label{eq:distal}
 	c_d(r) = \sigma^\prime   \int_{0}^{\infty}  \frac{J_0( \rho \,   r) J_0( \rho \,   L)\rho }{\rho^{2} + q}  \, d\rho 
 \end{equation} 
holds in the spatial domain. 
 In this case $ \sigma^\prime/\sigma= L I_1/I_2$, where \textit{I} denotes the value of the respective Bessel integrals at \textit{L}.
 The matching can be appreciated in Fig. \ref{fig:flapk3}. 
 
\begin{figure}[h]
 	\centering
 	
 	\includegraphics[width=1.0\linewidth]{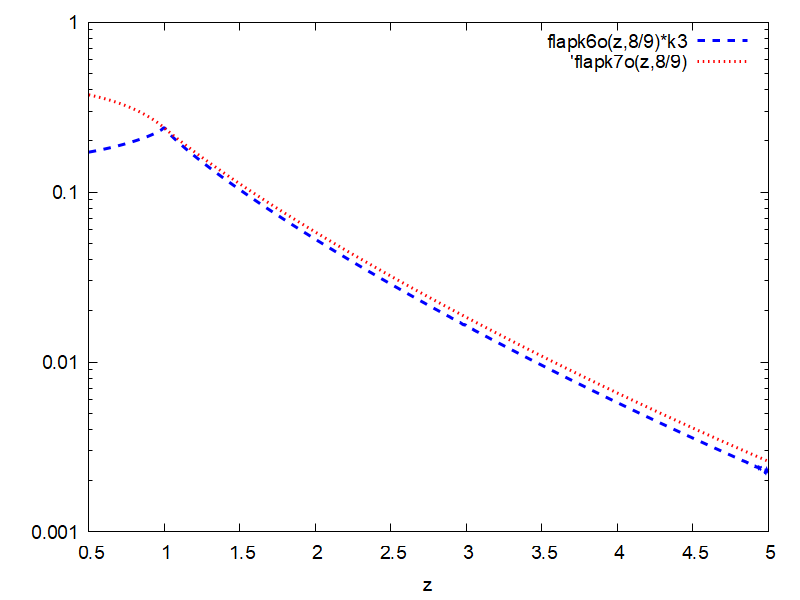}  
 	
 	The full solution (eq. \ref{eq:full}) is compared to the outer one (eq. \ref{eq:distal}) for $\alpha = 8/9$.
 	
 	\caption{Comparison of the fractional-order solutions}
 	\label{fig:flapk5}
 \end{figure}

 In the fractional-order case, due to the non-local character of the Riesz Laplacian this approximation does not recover the solution but achieves only an asymptotic similarity. 
 The differences in the obtained solutions can be appreciated in Fig. \ref{fig:flapk5}. 
 Eq. \ref{eq:distal} can be used as a definition of an asymptotic solution also in the fractional case where in the kernel denominator the exponent 2 is substituted by $2 \alpha$.
 
 A more tractable asymptotic solution can be obtained when one moves the impulse source to the origin -- i.e. $L=0$. In this case 
 \begin{equation}
 \hat{c} = \frac{\sigma }{|k|^{2\alpha } + q} 
 \end{equation}
 is obtained as a solution in the Fourier domain. 
 Therefore, the asymptotic solution is obtained as the Hankel transform:
 \begin{equation}\label{eq:dist}
 c_a(r) = \sigma   \int_{0}^{\infty}  \frac{J_0( \rho \,   r) \rho }{\rho^{a} + q}  \, d\rho, \quad a = 2 \alpha   
 \end{equation}
 The function is plotted in Fig. \ref{fig:flapk6} with $\alpha$ parametrization.

 Remarkably, for $a=2$  for the outer component we obtain   
 \begin{equation}\label{eq:besseljk}
 c_2(r)= \sigma^\prime   \int_{0}^{\infty}  \frac{J_0( \rho \,   r)  \rho }{\rho^{2 } + q} d\rho = \sigma^\prime K_0(\sqrt{q} r)
 \end{equation}
 where, $K_0$ is the corresponding modified Bessel function and $\sigma^\prime =k_2$ for this case. 
 In such way, the fractional problem leads to a special function, related to the Bessel $K_0$ function \citep{Prodanov2022}.

 \begin{figure}[h]
 	\centering
 	\includegraphics[width=1.0\linewidth]{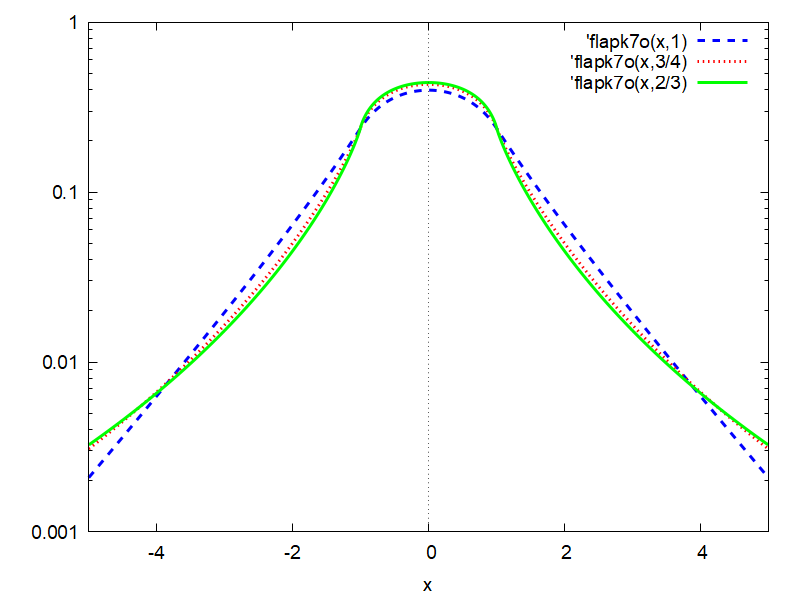}  
 	\caption{Influence of the fractional exponent on the shape of the full solution}	
 	\label{fig:flapk2}
 \end{figure}

	Furthermore, the asymptotic solution can be represented as Fox H-function by virtue of the Mellin transforms  theory (Appendix \ref{sec:fox}). 
	The Mellin transform of the asymptotic solution is given by the integral
	\begin{multline}
	C (s)=\mathcal{M}_x [c_a] (s) =\int_{0}^{\infty} x^{s-1}	c_d(x) dx = \\ \int_{0}^{\infty}  d\rho \int_{0}^{\infty} x^{s-1} \frac{J_0( \rho \,   x) \rho }{\rho^{a} + q}  \, d x 
	\end{multline}
	We use the known Mellin transform pair of the Bessel $J_0(x)$ function (see for example \cite{Oberhettinger1974})
	$$
	\mathcal{M}_x [J_0(\rho x)] (s) =
	\frac{2^{s-1}}{\rho^s} \frac{\Gamma \left(s/2 \right) }{\Gamma \left(1-s/2 \right) }
	$$
	to obtain
	\begin{multline}
	C (s)=\int_{0}^{\infty} \frac{\Gamma \left( \frac{s}{2}\right)  {2}^{s-1} {\rho}^{1-s}}
	{\Gamma \left( 1-\frac{s}{2}\right)  \left( {{\rho}^{a}}+q\right) } d \rho= \\ \frac{\Gamma \left( \frac{s}{2}\right)  {2}^{s-1}}{\Gamma \left( 1-\frac{s}{2}\right) } \int_{0}^{\infty} \frac{ {\rho}^{1-s}}
	{  {\rho}^{a}+q } d \rho  , 
	\end{multline}
	The last integral can be recognized as an Euler Beta  integral:
	\[
	\int_{0}^{\infty} \frac{ {\rho}^{1-s}}
	{  {\rho}^{a}+q } d \rho =  \frac{q^{(2-s)/a-1}}{a} \mathcal{B}\left(1- \frac{2-s}{a}, \frac{2-s}{a} \right) 
	\]
	Therefore, 
	\begin{multline}\label{eq:Cs}
		C (s)= \frac{ {{q}^{\frac{2-s}{a}-1}} \Gamma\left( 1-\frac{2-s}{a}\right)  \Gamma\left( \frac{2-s}{a}\right)  \Gamma\left( \frac{s}{2}\right)  {{2}^{s-1}}}{a  \Gamma\left( 1-\frac{s}{2}\right) }
	\end{multline}
	The above result allows one to  invert the expression by the Mellin-Barnes integral and the residue Theorem. 
	The s-independent pre-factor is given by $\lambda =\sqrt[a]{q^2}/a$.
	The equation defines a Fox H -function  with kernel
	\begin{equation}
	\mathcal{H}^{m,n}_{p,q}(s)=\frac{ \Gamma\left( \frac{s}{2}\right) \Gamma\left( \frac{2}{a} -\frac{s}{a}  \right) \Gamma\left( 1-\frac{2}{a} +\frac{s}{a}\right)     {{2}^{s-1}}}{  \Gamma\left( 1-\frac{s}{2}\right) }
	\end{equation}
	From where we read off the parameters $m=2$, $n=1$, $p=1$, $q=3$.
	Therefore, up to a constant factor $\lambda$, %
	the inverse Mellin transform is given by the H-function
	\begin{equation}\label{eq:foxasymp}
		c_a(r) =\lambda H^{2,1}_{1,3}\HF{\left(1-\frac{2}{a}, \frac{1}{a}\right)\Big|   }{ \left(0, \frac{1}{2} \right), \left(1-\frac{2}{a}, \frac{1}{a}\right) \Big| \left(0, \frac{1}{2} \right) }{\frac{\sqrt[a]{q} r}{2}}
	\end{equation}
	where $a = \beta+1$ and $\lambda = \frac{\sqrt[a]{q^2}}{ a q}$.
	For the given kernel one may take evaluation contour \textit{L} to be a Bromwich-type vertical line Re(s)=c with $-2< c < 2/a$ for $(a>1)$, indented to avoid poles, 
	which places the poles of $\Gamma(s/2)$ and $\Gamma(1-2/a+s/a)$  to the right and $\Gamma(2/a-s/a)$ poles to the left.

	It should be noted that for rational values of th exponent \textit{a} the H function is reducible to a Meijer G-function (see Appendix \ref{sec:mg}) and eventually represented by a finite combination of elementary and special functions.	This will be a subject of further studies. 
	For example, let a=2. Then two of the gamma factors cancel and one obtains the simpler kernel
	\begin{equation}
	\mathcal{H} = 2^{s-1} \Gamma^2 \left( \frac{s}{2}\right) 
	\end{equation}
	which is the Mellin transform pair of $K_0(z)$, thus confirming eq. \ref{eq:besseljk}. In G-function notation
		\begin{equation}
	  K_0(z) =\frac{1}{2} G^{2,0}_{0,2} \HF{-}{0,0}{\frac{z^2}{4}}
		\end{equation}

\begin{figure}
	\centering
	\includegraphics[width=1.0\linewidth]{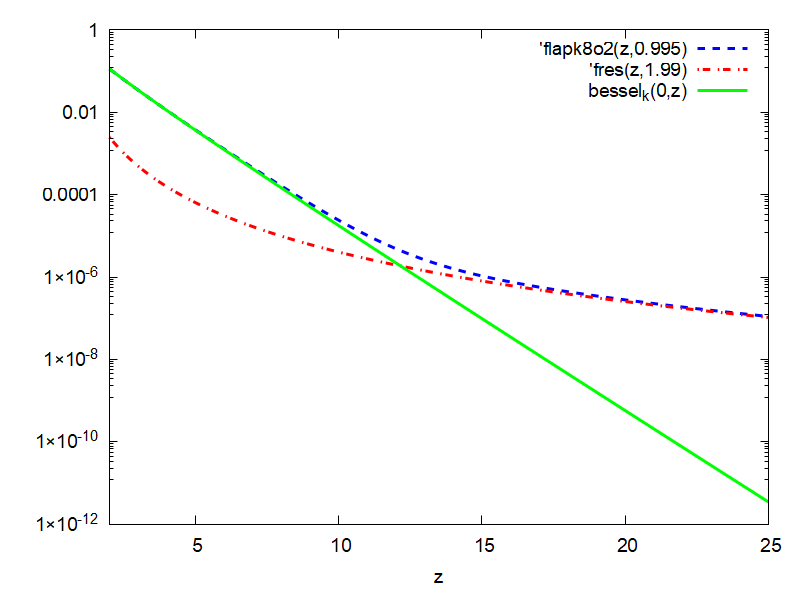}
	\caption{Asymptotic behavior of $c_a(z)$ for $\alpha=0.995$}
	\label{fig:flapk7}
\end{figure}

	Finally, so-identified H-function from eq. \ref{eq:foxasymp} allows one to derive  asymptotics for large values \textit{r} even without computing explicitly the function. 
	For the range of interest $a \in (1, 2]$ $c_a$ has a purely algebraic large-r expansion generated by the right–half‑plane poles
   located at $2-s= -a k, \ k \in \fclass{N}{}$. Therefore, the leading residue (i.e. for  $k=1$) is
	\begin{equation}\label{eq:asympfrac}
	c_1 = {\Gamma\left( \frac{a}{2} + 1\right) }^{2} \ \frac{  {{2}^{a+1}} \sin{\left( \frac{ \pi  a}{2}\right) } {}}{ \pi{z}^{ a+2} }, \quad a =2 \alpha 
	\end{equation}
 	A plot for $\alpha=0.995$ is presented in Fig. \ref{fig:flapk7}.
 
\section{Numerical experiments}\label{sec:numerics}
Numerical experiments were performed in the computer algebra system Maxima v. 5.47.0 running in a Windows (TM) 10 operation system. 
To this end, the Double-Exponential (DE) routine was ported to Lisp and integrated into the computer algebra system Maxima. 
The code can be downloaded from \url{https://github.com/dprodanov/intde}.
Plots of the solutions have been computed and rendered using the DE integration routine.  Everywhere $q=1$ was used. 

\subsection{The Double-Exponential quadrature method}\label{sec:de}
The DE quadrature integration  \citep{Takahasi1974,Mori1985} can be summarized as follows.
The integral on the interval $A$ 
\begin{equation}
	I= \int_{A} f(x) dx = \int_{A} f(x) dx= 
	\int_{-\infty}^{\infty} f \circ \phi \ (t) \phi^{\prime} (t) dt 
\end{equation}
is transformed to an integral on the entire real line
and the function $\phi$ guarantees double exponential convergence to  \textit{I} as the limit of the Riemannian sum
\begin{equation}
	I= 	\llim{N}{\infty}{ } S_N =
	\llim{N}{\infty}{ }  h \sum\limits_{k=-N}^{N} \underbrace{\phi^{\prime}(k h)}_{w_k} f(\underbrace{\phi(k h)}_{x_k})  
\end{equation}
Therefore, the value of the integral can be approximated by the truncated sum $S_N \approx I$, where the weights are given by $w_k$ and the abscissas (i.e. evaluation points) by $x_k$, while $h$ is an adjustable parameter. 
For the case of a semi-finite interval the DE method employs the transformation
\begin{equation}
	x  = \phi(t)= \exp{\left( \frac{\pi}{2}\sinh{t}\right) }
\end{equation}
The absolute error of the method is of the order $\exp{(–c\, N/\log{N})}$. 
Oscillatory kernels of the form $f(x) \sim g(x)  \sin(\omega   x + \theta)$ at infinity are supported by the modification of DE method \citep{Ooura1991} using
 \begin{equation}
 	x = \phi(t)= \frac{t}{1 -\exp{\left( -6 \sinh{t}\right) }} 
 \end{equation}

\begin{figure}
	\centering
	\includegraphics[width=1.0\linewidth]{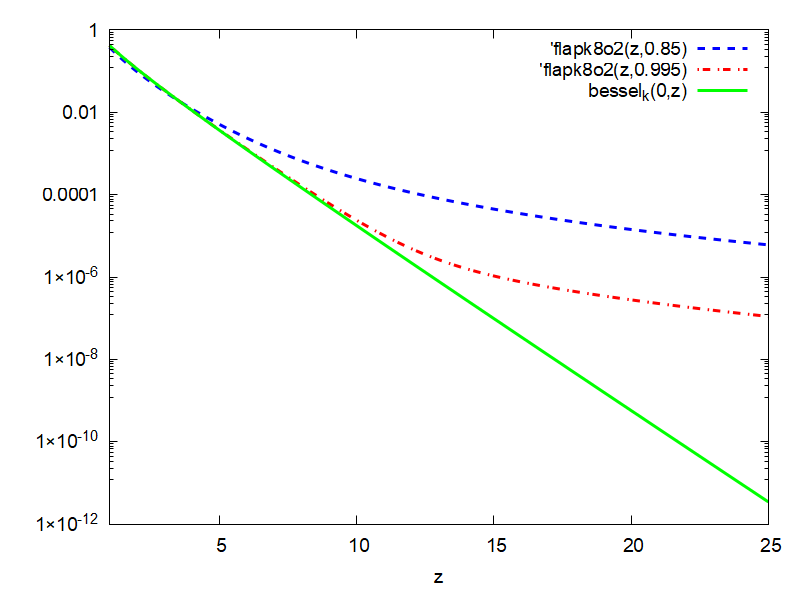}
	\caption{Comparison of $c_a(z)$ for $\alpha$=0.85, 0.995  with $K_0 (z)$ }
	\label{fig:flapk6}
\end{figure}
The outer asymptotic solution $c_a (r)$ (eq. \ref{eq:dist}) is compared for two values of $\alpha$ -- 0.85 and 0.995 -- with the integer-order one in Fig. \ref{fig:flapk6}.  

\subsection{Convergence acceleration method for numerical Hankel transforms}
A way to overcome the slow convergence of the Bessel integrals 
is to evaluate intermediate integrals between the zeros of the Bessel J function \citep{Lucas1995}  
 and to apply convergence acceleration techniques on the partial sums. 
 Even more economic approach is to use the approximation of the Bessel J zeros which can be computed as
\begin{equation}
r_{\nu} (k) = \pi\left( k+ \frac{\nu}{2}  - \frac{1}{4} \right)  + O\left(  {1}/{k} \right) 
\end{equation}
for large arguments; here $\nu$ denotes the order of the Bessel function.
An even improved estimate from the same reference is
 \begin{multline}
 	r_{\nu} (k) = b - (4 n^2 - 1 )/(8 b) + O(1/k^2),  \\ b=\pi\left( k+ \frac{\nu}{2}  - \frac{1}{4} \right)   
 \end{multline}
Then the Hankel transform of $F$ can be evaluated for some large integer \textit{N} as the sum
\begin{multline}
	f(z) = \sum_{k=0}^{N} \int_{r_\nu(k)}^{r_\nu(k+1)} F(r) J_{\nu} (r z) r dr + \\ \int_{r_\nu(N+1)}^{\infty} F(r) J_{\nu} (r z) r dr
\end{multline}
 where $r_\nu(0)=0$ is formally adjoined to the zero set. 
The last term of the sum can be used as an error estimate or, alternatively, it can be approximated using the Bessel function asymptotic.  
A method based on the DE transformation and employing Wynn's $\epsilon$-algorithm for convergence acceleration has been introduced in \cite{Prodanov2022}.  
An acceptable approximation was found to use only the first 15 asymptotic zeros.
The algorithm was employed to plot Fig. \ref{fig:flapk6}.

\subsection{Plots of the solutions}
Eqs. \ref{eq:full} and \ref{eq:distal}, representing the solutions,  can be computed numerically to a predefined order  of precision.
Plots of the solutions are presented in Figs.   \ref{fig:flapk3}-- \ref{fig:flapk2}  for $L=1$, $s=1$, and $q=1$ using the  DE method for oscillatory integrals.
The full solution is plotted for $\alpha =1, 3/4, 2/3$ in Fig. \ref{fig:flapk2}.
The impact of the fractional exponent can be appreciated in Figs. \ref{fig:flapk7} and \ref{fig:flapk6}.
From the plot  the heavy trails of the concentration can be appreciated against the baseline of $K_0$ (i.e. integer-order solution).

\section{Discussion}\label{sec:disc}

The fractional reaction-diffusion system (eq. \ref{eq:syst}) has been previously  used to model the distribution of  diffusing species around an implanted electrode  \citep{prodanov2016c}.
Such electrodes are routinely used in neurophysiological studies and deep brain stimulation  for sensing neural activity or deep brain stimulation for therapeutic applications. 
Typical implanted electrodes have circular or rectangular cross-section shapes  and high aspect ratios.

The fractional Laplacian can be defined on unbounded domains in  several equivalent ways  \citep{kwasnicki2017}. However, when these definitions
are restricted to bounded domains, the associated boundary conditions lead to distinct operators (discussion in \cite{Lischke2020}). 
Although reasoned through the spectral representation (Appendix \ref{sec:riesz}), 
the present paper   favors the use of a potential approach implementing long-range spatial interactions  since it connects with the physical interpretation of the diffusion equation.
Although the presented paper employs the Riesz Laplacian, the modeling approach is not restricted by the choice of operator and can be extended to other definitions.  

The present contribution obtained an analytical solution of the problem, which is computable by quadrature (eq. \ref{eq:full}) as demonstrated in Fig. \ref{fig:flapk2}.
From theoretical perspective, of some interest is also
the asymptotic solution (eq. \ref{eq:foxasymp}) obtained in terms of a Fox H-function. 
Computation of such functions is an open area of research since they are very general objects. 

Finally, of practical interest is also  eq.  \ref{eq:asympfrac} (see also Fig. \ref {fig:flapk7}).
Its application would allow for a more flexible approach when estimating transport parameters from experimental data.


%

\begin{ack}
The present work is funded by the Horizon Europe project VIBraTE, grant agreement no. 101086815. 
\end{ack}

\appendix
\section{Integral Transforms}
\subsection{Fourier transform}\label{sec:ft}
The Fourier transform will be defined under the "engineering" convention
\[
\hat{f}(k)= \mathcal{F} f(x) := \int_{\fclass{R}{d}} f(x) e^{-i k \cdot x} dx^d
\]
with an  inverse
\[
f(x) = \mathcal{F}^{-1} \hat{f}(k):=\frac{1}{(2 \pi)^d} \int_{\fclass{R}{d}} \hat{f}(k) e^{ i k \cdot x} dx^d
\]
\subsection{Hankel transform convention}\label{sec:ht}
The automorphic Hankel transform is defined as \citep{poularikas2010}:
\begin{equation}
\hat{f} (\rho)=\mathcal{H}_{\nu} f (z):= \int_{0}^{\infty}  f(z) J_{\nu}( \rho \,   z) z dz   
\end{equation}
Remarkably, for 2 spatial dimensions the Laplacian is represented by a monomial factor:
\[
\mathcal{H}_{0} \Delta f (z) = - \rho^2 \hat{f} (\rho)
\]
and in the similar way the image of the Riesz Laplacian is 
$
(-\Delta)^{\alpha} \mapsto - |\rho|^{2 \alpha}
$.

\subsection{Mellin transform}

The Mellin transform of a function $f(t)$ is defined as
\begin{equation}
\mathcal{M}_t[f](s) := \int_{0}^{\infty} t^{s-1} f(t) dt =F (s),
\end{equation}
wherever the integral exists on the complex plane.
The inverse Mellin transform is given by complex integration along a Bromwitch contour:
\[
f(z)=	\mathcal{M}^{-1}[F](z)= \frac{1}{2\pi i} \int_{Br}  F(s)  z^{-s} ds 
\]

\section{The Riesz Laplacian operator}\label{sec:riesz}
The Riesz operator can be defined in the Fourier domain \citep{kwasnicki2017}  where 
$$ - (-\Delta)^{\alpha} f(x) \mapsto   |k|^{2 \alpha} \hat{f}(k). $$
Then starting from the algebraic  substitution 
$
(-\Delta)^{\alpha} \mapsto - |k|^{2 \alpha}
$, the fractional Laplacian can be interpreted as the gradient of another operator, that is
\[
- |k|^{2 \alpha} = i \mathbf{k} \cdot i\mathbf{k^0} |k|^{2 \alpha-1} , \quad  \mathbf{k^0} = \frac{\mathbf{k}}{k}
\]
where the dot denotes the scalar product and  $\mathbf{k^0} $ is a unit wave vector  in the Fourier space.
Therefore, the associated Riesz gradient can be defined as a pseudo-differential operator  
\[
\nabla^\beta \mapsto  i\mathbf{k^0} |k|^{1-\beta} = i\mathbf{k} /|k|^{\beta}, \quad  \beta= 2 \alpha -1
\]
for a suitable function space. This has the advantageous physical interpretation of a \textit{generalized first Fick's law}, which 
 This corresponds to the convolution with 
 \[
 \nabla^\beta = I^{1-\beta} \ast \nabla
 \]
 where $ I^{1-\beta}$ denotes the Riesz potential.

\section{Fox H functions}\label{sec:fox}
The Fox H-function
$H^{m,n}_{p,q}\left[ z |\cdot \right]$  is defined by a Mellin-Barns integral as follows:
\begin{multline}
	H^{m,n}_{p,q} \left[ z \middle| \begin{array}{l} (a_1, A_1), (a_2, A_2), \ldots, (a_m, A_m) \\ (b_1, B_1), (b_2, B_2), \ldots, (b_n, B_n) \end{array} \right] := \\	\frac{1}{2\pi i} \int_L \mathcal{H}^{m,n}_{p,q}(s) \cdot z^{-s} ds =
	\frac{1}{2\pi i} \int_L \frac{\prod_{j=1}^m \Gamma(b_j + B_j s)}{\prod_{j=1+m}^q \Gamma(1-b_j - B_j s)} \cdot \\
	\frac{\prod_{j=1}^n \Gamma(1-a_j - A_j s)}{\prod_{j=n+1}^p \Gamma(a_j + A_j s)} \cdot z^{-s} ds,
\end{multline}
where:
\begin{itemize}
	\item 	$m$ and $n$ are non-negative integers,
	\item  $p$ and $q$ are positive integers,
	\item $z$ is a complex variable,
	\item $a_j, b_j, A_j, B_j$ are  parameters with positive real parts,
\end{itemize}
whereas $L$ is a suitable contour in the complex plane that separates the poles of the Gamma functions in the numerator.

\section{Meijer G functions}\label{sec:mg}
The Meijer G functions are closely related to the Fox H-functions.
The G function $G^{m,n}_{p,q}\left[ z |\cdot \right]$  is defined by a Mellin-Barns integral as follows:
\begin{multline*}
	G^{m,n}_{p,q} \left[ z \middle| \begin{array}{l} a_1, a_2, \ldots, a_p \\ b_1, b_2, \ldots, b_q \end{array} \right] := \\	\frac{1}{2\pi i} \int_L \mathcal{G}^{m,n}_{p,q}(s) \cdot z^{-s} ds =
	\frac{1}{2\pi i} \int_L \frac{\prod_{j=1}^m \Gamma(b_j - s)}{\prod_{j=1+m}^q \Gamma(1-b_j +s)} \cdot \\
	\frac{\prod_{j=1}^n \Gamma(1-a_j + s)}{\prod_{j=n+1}^p \Gamma(a_j -s)} \cdot z^{-s} ds,
\end{multline*}
where:
\begin{itemize}
	\item 	$m$ and $n$ are non-negative integers,
	\item  $p$ and $q$ are positive integers,
	\item $z$ is a complex variable,
	\item $a_j, b_j$ are  parameters with positive real parts, and
\end{itemize}
and the contour $L$ separates poles. 
\bibliography{mechdiffusion}             
\end{document}